\documentclass[11pt]{article}
\usepackage{graphicx}
\usepackage{amsmath}
\usepackage{amsfonts}
\usepackage{epsfig}
\usepackage{setspace}

\newcommand{\be}{\begin{equation}}
\newcommand{\ee}{\end{equation}}
\newcommand{\beqn}{\begin{eqnarray}}
\newcommand{\eeqn}{\end{eqnarray}}
\newcommand{\beqns}{\begin{eqnarray*}}
\newcommand{\eeqns}{\end{eqnarray*}}

\newcommand{\Var}{\mbox{Var}}

\newcommand{\Cov}{\mbox{Cov}\ }

\newcommand{\EE}{\ensuremath{{\mathbb E}}}
\newcommand{\II}{\ensuremath{{\mathbb I}}}

\newcommand{\fr}[1]{(\ref{#1})}

\newtheorem{lemma}{Lemma}
\newtheorem{theorem}{Theorem}

\newtheorem{remark}{Remark}

\begin{document}

\title{\Large{\bf Minimax adaptive wavelet estimator for the simultaneous blind deconvolution with fractional Gaussian noise}}

\author{
\large{ Rida Benhaddou}  \footnote{E-mail address: Benhaddo@ohio.edu}
  \\ \\
Department of Mathematics, Ohio University, Athens, OH 45701} 
\date{}

\doublespacing
\maketitle
\begin{abstract}
We construct an adaptive wavelet estimator that attains minimax near-optimal rates in a wide range of Besov balls. The convergence rates are affected only by the weakest dependence amongst the channels, and take into account both noise sources. \\


{\bf AMS (2000) Subject Classification: 62G05, 62G20, 62G08 }
 \end{abstract} 

\section{Introduction.}

Consider the problem of estimating the unknown response function $f$ based on the noisy convolutions that are continuously observed as $Y_l(t)$, $l=1, 2, \cdots, M$, and described by the model 
\be \label{conveq}
Y_l(t)=\int^1_0f(s)g_l(t-s)ds + \varepsilon^{\alpha_{1l}} Z_1^{H_{1l}}(t), \ \  t \in [0, 1],
\ee
where $f(t)$ is periodic, $Z_1^{H_{1l}}(t)$ are independent fractional Gaussian processes,  $\alpha_{1l}=2-2H_{1l}\in (0, 1]$ are the parameters of long-range-dependence (LRD), and $H_{1l}$ are Hurst parameters, $l=1, 2, \cdots, M$. The kernel functions $g_l(t)$ are unknown.  Instead, one continuously observes 
\be \
g_l^{\delta}(t)=g_l(t) + \delta^{\alpha_{2l}} Z_2^{H_{2l}}(t), \ \  t \in [0, 1], l=1, 2, \cdots, M, \label{kernel}
\ee
where $Z_2^{H_{2l}}(t)$ are independent fractional Gaussian processes, with $\alpha_{2l}=2-2H_{2l}\in (0, 1]$. The quantities $Z_1^{H_{1l}}(t)$ and $Z_2^{H_{2l}}(t)$ are assumed to be independent of each other.  The objective is to estimate $f(t)$. This is another version of blind deconvolution. 

Inverse problems  with unknown operators in their general aspect have been studied by Hoffmann and Reiss~(2008) where two different approaches are suggested to handle the issue. Delattre et al.~(2012) implement a blockwise Singular Value Decomposition (SVP) to treat the blind deconvolution problem when the signal belongs to some Sobolev class. Vareschi~(2015) looks into the Laplace deconvolution with noisy kernel when the signal belongs to a Laguerre-Sobolev class. Benhaddou~(2018a) investigated the blind deconvolution model with fractional Gaussian noise (fGn) for the one channel standard deconvolution when the signal is periodic and belongs to some Besov class and the kernel is contaminated with white noise. Recently, Benhaddou~(2018b) derived the lower bounds for the wavelet estimators under the exact same setting as in the present work.

Standard (Fourier) deconvolution model has witnessed the publication of a great deal of papers, including those that deal with the issue of long-memory (LM) or long-range dependence (LRD). A detailed literature review on that can be found in Benhaddou~(2018b).  

The case $\delta=0$ and $\alpha_{1l}=1$, $l=1, 2, \cdots, M$, corresponds to the simultaneous deconvolution with white noise and known kernel studied in De Canditiis and Pensky~(2006), while the case $\delta=0$ corresponds to the multichannel deconvolution model with long-range dependence and known kernel investigated in Kulik et al.~(2015). In addition, the case $M=1$ and $\alpha_{2l}=1$, pertains to the one channel blind deconvolution model with fGn investigated in Benhaddou~(2018a), while the case $\delta=0$ and $M=1$, corresponds to the one channel deconvolution model with fBm and known kernel discussed in Wishart~(2013).  

The objective of the paper is to complement the work in Benhaddou~(2018b) by constructing an adaptive hard-thresholding wavelet estimator for model \fr{conveq}. We focus on the regular-smooth convolution and, following Wishart~(2013), we apply Wavelet-Vaguelette-Decomposition (WVD) via Meyer-type wavelets to de-correlate fGn. In addition, similar to Benhaddou~(2018a), a preliminary stabilizing  thresholding procedure is applied to estimate the wavelet coefficients, and the standard hard-thresholds are then applied to keep only the coefficients of the wavelet expansion so as the variance is minimal. We show that the proposed approach is asymptotically near-optimal over a wide range of Besov balls under the $L^2$-risk. In addition, we show that the convergence rates are expressed as the maxima between two terms, taking into account both the noise sources. Moreover, the convergence rates depend only on the largest long-memory parameters, $\alpha_{il}$, $l=1, 2, \cdots, M$, which correspond to the weakest dependence from amongst the $M$ channels. These rates deteriorate as $\max_{l\leq M}\{\alpha_{il}, i=1,2\}$ get smaller and smaller. Similar behavior has been pointed out in Wang~(1997), Wishart~(2013), Kulik et al.~(2015) and Benhaddou~(2018a). It should be noted that with $\delta=0$, our convergence rates are similar to those in Kulik et al.~(2015), and with $\delta=0$ and $\alpha_{il}=1$, $i=1,2$, $l=1, 2, \cdots, M$, our convergence rates match those in De Canditiis and Pensky~(2006). Finally, with $M=1$ and $\alpha_{2l}=1$, our rates coincide with those in Benhaddou~(2018a), while with $\delta=0$ and $M=1$, our rates match up so some logarithmic factor of $\varepsilon\asymp n^{-1/2}$, those in Wishart~(2013).
 \section{Estimation Algorithm.}
  In what follows, denote $U=[0, 1]$, and
 let  $\tilde{h}(m)$ be Fourier coefficient of the function $h(t)$. Also, let $a\vee b = \max(a,b)$ and $a \wedge b = \min(a,b)$. Consider a Meyer-type wavelet basis $\psi_{j, k}(t)$ and let $m_{0}$ be its lowest resolution level and denote the scaling function for the wavelet by $\psi_{m_0-1, k}(t)$. Since the functions $\psi_{j, k}(t)$ form orthonormal bases of the $L^2(U)$ space, the function $f(t)$ can be expanded over these bases with coefficients $\beta_{j, k}$ into wavelet series as
 \be
f(t)= \sum^{\infty}_{j=m_0-1} \sum^{2^j-1}_{k=0}\beta_{j, k}\psi_{j, k}(t).
\ee
Applying Fourier transform to equations \fr{conveq} and \fr{kernel} yields 
\beqn
\tilde{Y}_l(m)&=& \tilde{f}(m)\tilde{g_l}(m)+ \varepsilon^{\alpha_{1l}}\tilde{Z}^{H_{1l}}_1(m),\label{ytild}\\
\tilde{g}^{\delta}_l(m)&=&\tilde{g}_l(m) + \delta^{\alpha_{2l}}\tilde{Z}^{H_{2l}}_2(m).\label{gtild}
\eeqn
For the Fourier coefficients of $f(t)$, $\tilde{f}(m)$, consider the weighted estimators given by 
\be \label{g-trunc}
\widehat{\tilde{f}(m)}= \left\{ \begin{array}{ll} {\frac{\sum^M_{l=1}\omega_{l}(m)\tilde{g}^{\delta}_l(m)\tilde{Y}_l(m)}{\sum^M_{l=1}\omega_{l}(m)|\tilde{g}^{\delta}_l(m)|^2}}, & \mbox{if}\ \ \min_{l \leq M}|\tilde{g}^{\delta}_l(m)|^2 >{k^2} \delta^{2\alpha^*_{2}}|m|^{\alpha^*_{2}-1} |\ln(\delta)|,\\
0,&  \mbox{if}\ \ otherwise,
\end{array} \right.
\ee
where $k$ is a positive constant independent of $m$ and $\delta$, $\omega_{l}(m)$ are weights to be determined later, and $\alpha^*_2=\max\{\alpha_{21}, \alpha_{22}, \cdots, \alpha_{2M}\}$. If $\psi_{j, k, m}=<e_m, \psi_{j, k}>$ are Fourier coefficients of $\psi_{j, k}(t)$, then, by Plancherel formula and \fr{g-trunc}, we obtain the truncated estimator
\be \label{be-estim}
\tilde{\beta}_{j, k}= \sum_{m\in W_j} \widehat{\tilde{f}(m)}\psi_{j, k, m},
\ee
where, for any $j \geq m_0$, 
\be \label{omeg}
W_j=\left \{ m: \psi_{j, k,m}\neq 0 \right\} \subseteq 2\pi/3\left[ -2^{j+2}, -2^{j} \right] \cup \left[ 2^j, 2^{j+2} \right ],
\ee
since Meyer wavelets are band-limited (see, e.g., Johnstone et al.~(2004)).  Then, define the estimator for $f(t)$ as
\be \label{ef-hat}
\widehat{f}_{\varepsilon, \delta}(t)= \sum^{J-1}_{j=m_0-1} \sum^{2^j-1}_{k=0} \widehat{\beta}_{j, k}\psi_{j, k}(t),
\ee
where
\be \label{be-hat} 
 \widehat{\beta}_{j, k}= \tilde{\beta}_{j, k} \II \left(|\tilde{\beta}_{j, k}| > \lambda^{\alpha}_{j; \varepsilon, \delta} \right),
\ee
and the values of $J$, $m_0$ and $\lambda^{\alpha}_{j; \varepsilon, \delta}$ are to be determined.  Next we introduce a condition that the functions $g_l(t)$ satisfy.\\  
\noindent
{\bf Assumption 1.} The Fourier coefficients  $\tilde{g}_l(m)$ of kernels $g_l(t)$ are such that 
\be \label{kern}
c_{l1}|m|^{-2\nu_l} < |\tilde{g}_l(m)|^2 < c_{l2}|m|^{-2\nu_l}, \ \ \ l=1, 2, 3, \cdots, M,
\ee 
where $\nu_l>0$, $c_{l1}$ and $c_{l2}$ are some positive constants independent of $m$.\\
To determine the choices of  $J$, $m_0$ and $ \lambda^{\alpha}_{j;\varepsilon, \delta}$ in \fr{ef-hat} and \fr{be-hat}, it is necessary to evaluate the variance of \fr{be-estim}. Thus, recall that by \fr{omeg}, one has $|m| \asymp 2^j$, and define for some constant $0<\rho < 1/2$, the sets $\Omega_1$ and $\Omega_2$ as 
\beqn 
\Omega_1&=& \left \{ m \in W_j: \min_{l \leq M}|\tilde{g}^{\delta}_l(m)|^2 >{k^2} \delta^{2\alpha^*_{2}}|m|^{\alpha^*_{2}-1} |\ln(\delta)| \right\}, \label{Neumn1}\\ 
\Omega_2&= &\left \{ m \in W_j: \max_{l\leq M}|\delta^{\alpha_{2l}} Z^{H_{2l}}(m)|^2 < {\rho^2 k^2}\delta^{2\alpha^*_{2}} |m|^{\alpha^*_{2}-1}|\ln(\delta)| \right\}. \label{Neumn2}
\eeqn 
 Denote $\Omega_j=\Omega_1\cap\Omega_2$, and notice that on $\Omega_j$ one has
 \be  \label{blur-del}
\frac{ 1-2\rho}{1-\rho}|\tilde{g}_l(m)| \leq |\tilde{g}^{\delta}_l(m)| \leq \frac{1}{1-\rho} |\tilde{g}_l(m)|.
\ee
The next statement holds. 
\begin{lemma} \label{lem:Var}
Let $\tilde{ \beta}_{j, k}$ be defined in \fr{be-estim}. Choose the weights $\omega_{l}(m)$ in \fr{g-trunc} as 
\be \label{wei}
\omega_{l}(m)= \left(\varepsilon^{2\alpha_{1l}}|m|^{\alpha_{1l}-1}+ \delta^{2\alpha_{2l}}|m|^{\alpha_{2l}-1}\right)^{-1}.
\ee
Then, on $\Omega_j$ and under condition \fr{kern}, one has 
\be \label{var}
\EE\left| \tilde{ \beta}_{j, k}- \beta_{j, k}\right|^2 \asymp \frac{1}{M}\sum^M_{l=1}\left(\varepsilon^{2\alpha_{1l}}2^{j(2\nu_l+\alpha_{1l}-1)} + \delta^{2\alpha_{2l}}2^{j(2\nu_l+\alpha_{2l}-1)} \right),
\ee
and 
\be \label{var2}
\EE\left| \tilde{ \beta}_{j, k}- \beta_{j, k}\right|^4 \asymp \frac{1}{M^2}\sum^M_{l=1}\left(\varepsilon^{4\alpha_{1l}}2^{2j(2\nu_l+\alpha_{1l}-1)} + \delta^{4\alpha_{2l}}2^{2j(2\nu_l+\alpha_{2l}-1)} \right). 
\ee
\end{lemma}
Since the degrees of ill-posedness $\nu_1, \nu_2, \cdots, \nu_M$, are unknown, data-driven thresholds $ \lambda^{\alpha}_{j;\varepsilon, \delta}$ are necessary to make the estimator \fr{ef-hat} adaptive. Therefore, define the quantities
\be \label{sgdm}
S_j\left( \tilde{g}_l(m)\right)=  \sum_{m\in \Omega_j}|\tilde{g}_l(m)|^2,
\ee 
and notice that $2^j\left[S_j\left(\tilde{g}^{\delta}_l(m)\right)\right]^{-1}\asymp 2^{2j\nu_l}$. Following Lemma \ref{lem:Var} and \fr{sgdm} we choose the thresholds $\lambda^{\alpha}_{j;\varepsilon, \delta}$  of the form
\be \label{Thres2}
\lambda^{\alpha}_{j;\varepsilon, \delta} =
  \rho_1  \left[S_j\left(\tilde{g}^{\delta}_{l^*_1}(m)\right)\right]^{-\frac{1}{2}} \varepsilon^{\alpha_{1l^*_1}} {|\ln(\varepsilon)|^{\frac{1}{2}}}2^{{j\alpha_{1l^*_1}}/{2}}\vee \rho_2 \left[S_j\left(\tilde{g}^{\delta}_{l^*_2}(m)\right)\right]^{-\frac{1}{2}} \delta^{\alpha_{2l^*_2}} |\ln(\delta)|2^{j\alpha_{2l^*_2}/2}, 
 \ee
 where
 \beqn \label{oplsurr1}
l^*_{1}&=& \arg\min_{1 \leq l\leq M}\left\{ \varepsilon^{2\alpha_{1l}}2^{j(\alpha_{1l} +1)}\left[S_j\left(\tilde{g}^{\delta}_l(m)\right)\right]^{-1}\right\},\\
l^*_{2}&=& \arg\min_{1 \leq l\leq M}\left\{ \delta^{2\alpha_{2l}}2^{j(\alpha_{2l} +1)}\left[S_j\left(\tilde{g}^{\delta}_l(m)\right)\right]^{-1}\right\},
\eeqn
 for any $j \geq m_0$. Based on \fr{Thres2}, choose $m_0$  and $J$ such that 
\be  \label{Lev:J}
2^{m_0}= |\ln(\varepsilon)|\wedge|\ln(\delta)|, \ \ 2^J=2^{J_1}\wedge 2^{J_2}, 
\ee
where 
\be  \label{Lev:Jad}
 {J_i}=  \max\left\{  j : \left[S_j\left(\tilde{g}^{\delta}_{l^*_i}(m)\right)2^{-j(\alpha_{il^*_1} +1)}\right]^{-1}\leq \Gamma^i_{\varepsilon, \delta}\right\}, i=1, 2,
\ee
and 
\be
\Gamma^1_{\varepsilon, \delta}= \left[ \frac{\varepsilon^{2\alpha_{1l^*_1}}}{A^2M}\right]^{-1}, \ \  \Gamma^2_{\varepsilon, \delta}=\left[ \frac{\delta^{2\alpha_{2l^*_2}}}{A^2M}\right]^{-1}.
\ee
Remark that by \fr{Lev:J} and \fr{Lev:Jad},  $J$ satisfies
\be  \label{Lev:Jasym}
 2^{J} \asymp \left\{ \left[ \frac{\varepsilon^{2\alpha_{1l^*_1}}}{A^2M}\right]^{-\frac{1}{2\nu_{l^*_1} + \alpha_{1l^*_1}}} \wedge\left[ \frac{\delta^{2\alpha_{2l^*_2}}}{A^2M}\right]^{-\frac{1}{2\nu_{l^*_2} + \alpha_{2l^*_2}}}\right\},
\ee
and $l^*_1$ and $l^*_2$ are such that 
\beqn \label{oplsurr2}
l^*_{1}&=& \arg\min_{1 \leq l\leq M}\left\{ \varepsilon^{2\alpha_{1l}}2^{(\alpha_{1l} +2\nu_l)}\right\},\\
l^*_{2}&=& \arg\min_{1 \leq l\leq M}\left\{ \delta^{2\alpha_{2l}}2^{(\alpha_{2l} +2\nu_l)}\right\}.
\eeqn
\section{Minimax adaptivity and convergence rates in the $L^2$-risk.}
\noindent
{\bf Assumption 2.}  Denote $s^*=s + 1/2 -1/p$, and assume that $f(t)$ belongs to the one-dimensional Besov ball; that is, its wavelet coefficients satisfy
\be  \label{Beso}
 B^{s}_{p, q}(A)=\left \{ f \in L^2(U): \left( \sum_{j} 2^{js^{*}q}\left (\sum_{k}| \beta_{j, k}|^p\right)^{q/p}\right )^{1/q} \leq A\right \}.
\ee
 It remains to see how estimator  \fr{ef-hat} performs in the minimax sense, so we evaluate the minimax convergence rates of \fr{ef-hat} for the $L^2$-risk. Define such risk over the set $\Theta$ as 
\be  \label{Besbl}
R^{2}_{\varepsilon}(\Theta)= \inf_{\tilde{f}} \sup_{f\in \Theta}\EE\|\tilde{f}_{\varepsilon}-f\|^{2}_{2},
\ee
where the infimum is taken over all possible estimators $\tilde{f}$ of $f$. The derivation of upper bounds of the $L^2$-risk relies on the following lemma.
\begin{lemma} \label{lem:Lar-D} 
Let $\tilde{ \beta}_{j, k}$ and $\lambda^{\alpha}_{\varepsilon, \delta}$ be defined by \fr{be-estim} and \fr{Thres2}, respectively. Define, for some positive constant $\eta$, the set 
\be  \label{Set:Thet}
\Theta_{j, k, \eta}=\left \{ \Theta: |\tilde{ \beta}_{j, k}-\beta_{j, k}| > \eta\lambda^{\alpha}_{\varepsilon, \delta} \right\}.
\ee
 Then, on $\Omega_1\cap\Omega_2$ and under condition \fr{kern}, as $\varepsilon, \delta \rightarrow 0$, simultaneously, one has 
\be  \label{Largdev}
\Pr \left( \Theta_{j, k, \eta}\right)= O \left ( \left[\varepsilon^{2\alpha_{1l^*_1}}\right]^{\frac{\rho_1^2\eta^2}{32\alpha \sigma_{o1}^2}}\vee \left[\delta^{2\alpha_{2l^*_2}}\right]^{\frac{\rho_2^2\eta^2}{32 \sigma_{o2}^2}} \right),
\ee
where $\rho_1$ and $\rho_2 $ appear in \fr{Thres2}, and  $\sigma^2_{oi}= \frac{\alpha_{il^*_i}}{{c_{1l^*_i}}} \left( \frac{8\pi}{3}\right)^{(2\nu_{l^*_i} + \alpha_{il^*_i} -1)}$, $i=1,2$.
\end{lemma}
Then, the following statement is true. 
\begin{theorem} \label{th:upperbds}
Let $\widehat{f}(t)$ be the wavelet estimator in \fr{ef-hat}, with $J$ given by \fr{Lev:J}-\fr{Lev:Jad} and $\lambda^{\alpha}_{j;\varepsilon, \delta}$ given by \fr{Thres2}. Let $s \geq \max\{\frac{1}{p}, \frac{1}{2} \}$ with $1 \leq p,q \leq \infty$, and let conditions  \fr{kern} and \fr{Beso} hold. If $\rho_1$ and $\rho_2 $ in \fr{Thres2} are large enough, then, as $\varepsilon, \delta \rightarrow 0$, simultaneously, one has
 \be \label{upprbds}
R^2_{\varepsilon}( B^{s}_{p, q}(A)) \leq C A^{2}\left\{ \begin{array}{ll} 
 \left[\frac{ \varepsilon^{2\alpha_{1l^*_{1}}}|\ln(\varepsilon)|}{A^2M} \right]^{\frac{2s}{2s+2\nu_{l^*_{1}}+ \alpha_{1l^*_{1}}}} \vee   \left[\frac{ \delta^{2\alpha_{2l^*_{2}}}\ln^2(\delta)}{A^2M} \right]^{\frac{2s}{2s+2\nu_{l^*_{2}}+ \alpha_{2l^*_{2}}}}, & \mbox{if} \   \ s > s_1\vee s_2,\\
 \left[ \frac{ \varepsilon^{2\alpha_{1l^*_{1}}}|\ln(\varepsilon)|}{A^{2}M}\right]^{\frac{2s^*}{2s^* +2\nu_{l^*_{1}} +\alpha_{1l^*_{1}} -1}}\xi_1 \vee  \left[ \frac{ \delta^{2\alpha_{2l^*_{2}}}
\ln^2(\delta)}{A^{2}M}\right]^{\frac{2s^*}{2s^* +2\nu_{l^*_{o}} +\alpha_{2l^*_{2}} -1}}\xi_2, & \mbox{if} \   \ s\leq s_1\wedge s_2,\\
 \left[\frac{ \varepsilon^{2\alpha_{1l^*_{1}}}|\ln(\varepsilon)|}{A^2M} \right]^{\frac{2s}{2s+2\nu_{l^*_{1}}+ \alpha_{1l^*_{1}}}} \vee  \left[ \frac{ \delta^{2\alpha_{2l^*_{2}}}
\ln^2(\delta)}{A^{2}M}\right]^{\frac{2s^*}{2s^* +2\nu_{l^*_{2}} +\alpha_{2l^*_{2}} -1}}\xi_2, & \mbox{if} \     \ s_1< s\leq s_2,\\
 \left[ \frac{ \varepsilon^{2\alpha_{1l^*_{1}}}|\ln(\varepsilon)|}{A^{2}M}\right]^{\frac{2s^*}{2s^* +2\nu_{l^*_{1}} +\alpha_{1l^*_{1}} -1}}\xi_1 \vee   \left[\frac{ \delta^{2\alpha_{2l^*_{2}}}\ln^2(\delta)}{A^2M} \right]^{\frac{2s}{2s+2\nu_{l^*_{2}}+ \alpha_{2l^*_{2}}}}, & \mbox{if}\    \ s_2<s\leq s_1,
 \end{array} \right.
\ee
where $\xi_1$ and $\xi_2$  are defined as 
\be  
\xi_1 =\left[|\ln(\varepsilon)|\right]^{\II \left( s= s_1 \right)}, \quad \mbox{} \quad \xi_2 =\left[|\ln(\delta)|\right]^{\II \left( s= s_2 \right)},
\ee
and
\beqn \label{opl}
s_i&=&\left(\frac{1}{p}-\frac{1}{2}\right)\left(2\nu_{l^*_{i}} + \alpha_{il^*_{i}}\right), \ \ \ i=1, 2.
\eeqn
\end{theorem}

 \begin{remark}	
{\rm{
{\noindent {\bf{(i)}}\,
  The upper-bounds \fr{upprbds} match, up to some logarithmic factors of $\varepsilon$ or $\delta$, the lower-bounds derived in Benhaddou~(2018b), and therefore estimator \fr{ef-hat} is asymptotically near-optimal over a wide range of Besov balls $B^s_{p, q}(A)$. \\
\bf{(ii)}}\,
Our convergence rates are expressed as the maxima between two terms, taking into account both noise sources (the signals and the kernels). This behavior was { pointed out} in Hoffmann and Reiss~(2008), Vareschi~(2015),  Benhaddou~(2018a) and Benhaddou~(2018b). In addition, the convergence rates depend on the largest amongst the long-memory parameters $\alpha_{il}$, $l=1, \cdots, M$, $i=1,2$, which correspond to the weakest LRD amongst the $M$ available channels, and deteriorate as $\max_{l\leq M}\{\alpha_{il}, i=1,2\}$ get closer and closer to zero.\\
 \noindent {\bf{(iii)}}\, For $\delta=0$, our rates coincide, up to some logarithmic factor of $\varepsilon\asymp n^{-1/2}$, with the upper bounds obtained in Kulik et al.~(2015) in the regular-smooth convolution case.  \\
 \noindent {\bf{(iv)}}\,
  For $M=1$ and $\delta=0$, our rates coincide with those in Wishart~(2013), up to some logarithmic factor of $\varepsilon\asymp n^{-1/2}$. \\
  \noindent {\bf{(v)}}\, 
 For  $\alpha_{11}=\alpha_{12}=\cdots=\alpha_{1M}=1$ and $\delta=0$,  our rates match, up to some logarithmic factor of $\varepsilon\asymp n^{-1/2}$, with the upper bounds obtained in De Canditiis and Pensky~(2006) in their regular-smooth convolution case.\\
 \noindent {\bf{(vi)}}\,  For $M=1$ and $\alpha_{2l}=1$, our rates match exactly those in Benhaddou~(2018a).  \\
 \noindent {\bf{(vii)}}\, 
Note that in practice, for the proposed estimation algorithm to be computationally possible, the data   $g^{\delta}_l(t)$ and $Y_l(t)$ must be of equal sizes. Therefore we cannot claim that one will achieve the same convergence rates as if  $g_l(t)$ were known if data $g^{\delta}_l(t)$ are chosen to have relatively a larger size than data $Y_l(t)$, as it was previously suggested in Benhaddou~(2018b).\\
  \noindent {\bf{(viii)}}\, 
The choices of $J$ and $\lambda^{\alpha}_{j;\varepsilon, \delta}$ in \fr{Thres2} and \fr{Lev:Jad} are independent of the parameters of the Besov ball and the smoothness parameters, $\nu_l$ of the unknown kernels $g_l$, and therefore estimator \fr{ef-hat} is adaptive with respect to those parameters.\\
{\bf{(ix)}}\,
 Finally, note that the long-memory parameters $\alpha_{il}$, $l=1, \cdots, M$, $i=1,2$ may not be known in advance, but can be estimated from the data. There are quite a few methods that have been developed to estimate the parameter $\alpha$ for various forms of LRD, including fGn. A comprehensive list can be found in Taqqu et al.~(1995), Fischer and Akay~(1996), Pilgram and Kaplan~(1998) and Heath and Vivero~(2012). One strategy is to have $2n$ observations from models \fr{conveq} and \fr{kernel} and for each channel, use the first $n$ observations to estimate $\alpha$ via any of the methods available and then use the remaining $n$ observations to estimate $f$ with $\alpha_{il}$ replaced by their sampling counterparts. }}
\end{remark}
\section{Proofs. }
{\bf Proof of Lemma \ref{lem:Var}.}
Note that by conditioning on  $\Omega_j$, the variance of \fr{be-estim} is
\beqn  \label{Var1}
\EE\left| \tilde{ \beta}_{j, k}-\beta_{j, k}\right|^2
 &\leq&  2\left(\frac{1-\rho}{1-2\rho}\right)^2\EE\left[\sum_{m\in W_j}\overline{\psi_{j, k, m}}\frac{\sum^M_{l=1}\omega_l(m)\varepsilon^{\alpha_{1l}}\tilde{g}_l(m)\tilde{Z}^{H_{1l}}(m)}{\sum^M_{l=1}\omega_l(m)|\tilde{g}_l(m)|^2} \right]^2 \nonumber\\
 &+& 2 \left(\frac{1-\rho}{1-2\rho}\right)^2\EE\left[\sum_{m\in W_j}\overline{\psi_{j, k, m}}\frac{\sum^M_{l=1}\omega_l(m)\delta^{\alpha_{2l}}\tilde{g}_l(m)\tilde{f}(m)\tilde{Z}^{H_{2l}}(m)}{\sum^M_{l=1}\omega_l(m)|\tilde{g}_l(m)|^2}\right]^2.
\eeqn
To evaluate \fr{Var1}, we use the following result from Benhaddou~(2016)
\be
\left|\Cov(\tilde{Z}^{H_{il}}(m),\tilde{Z}^{H_{il}}(m') )\right|^2\leq 2\left|mm'\right|^{1-2{H_{il}}}, i=1,2.\label{Cov:bd}
\ee
Now, plugging \fr{Cov:bd} in \fr{Var1}, taking into account \fr{blur-del}, $|\psi_{j, k, m}| \leq 2^{-j/2}$ and the fact that $|\tilde{f}(m)|\leq 1$, yields 
\be  \label{Var2}
\EE\left| \tilde{ \beta}_{j, k}-\beta_{j, k}\right|^2\leq  C\sum_{m\in W_j}|{\psi_{j, k, m}}|^2\frac{\sum^M_{l=1}\omega^2_l(m)|\tilde{g}_l(m)|^2\left[\varepsilon^{\alpha_{1l}}|m|^{\alpha_{1l}-1}+\delta^{\alpha_{2l}}|m|^{\alpha_{2l}-1}\right]}{\left[\sum^M_{l=1}\omega_l(m)|\tilde{g}_l(m)|^2\right]^2}.
 \ee
 Finally, minimizing \fr{Var2}  with respect to the weights $\omega_l(m)$ yields \fr{wei}. Consequently, using condition \fr{kern} completes the proof of \fr{var}. 
 \\
 To prove \fr{var2}, note that conditional on $\Omega_j$, the quantities in square brackets of \fr{Var1} are centered Gaussian random variables. Hence,  using some properties of Gaussian, \fr{var2} follows.  
 $\Box$\\
{\bf The proof of Lemma \ref{lem:Lar-D}.} We use the same conditioning argument on $\Omega_j$, and recall the set $\Theta_{j, k, \gamma}$ defined in \fr{Set:Thet}. Then, 
\be  \label{prob}
\Pr \left( \Theta_{j, k, \gamma}\right)\leq P_1+ P_2,
\ee
with
\beqns
 P_1&=& \Pr \left(\left|\eta^{\varepsilon}_1\right|> \frac{\gamma}{2}\lambda^{\alpha}_{\varepsilon, \delta}\right),\nonumber\\
 P_2&=& \Pr \left(\left|\eta^{\delta}_2\right|> \frac{\gamma}{2}\lambda^{\alpha}_{\varepsilon, \delta}\right),\nonumber\\
   \eeqns
and $\eta^{\varepsilon}_1$ and $\eta^{\delta}_2$ are centered Gaussian random variables having variances of the orders 
\be \label{nvarn1}
\Var (\eta^{\varepsilon}_1)\asymp \sum_{j\in W_j}|\psi_{j, k, m}|^2\frac{\sum^M_{l=1}\omega^2_l(m)\sigma^2_{1l}\varepsilon^{2\alpha_{1l}}|\tilde{g}_l(m)|^2|m|^{\alpha_{1l}-1}}{\left[\sum^M_{l=1}\omega_l(m)|\tilde{g}_l(m)|^2\right]^2},
\ee
and 
\be \label{nvarn2}
\Var (\eta^{\delta}_2)\asymp \sum_{j\in W_j}|\psi_{j, k, m}|^2\frac{\sum^M_{l=1}\omega^2_l(m)\sigma^2_{2l}\delta^{2\alpha_{2l}}|\tilde{g}_l(m)|^2|m|^{\alpha_{2l}-1}}{\left[\sum^M_{l=1}\omega_l(m)|\tilde{g}_l(m)|^2\right]^2}.
\ee 
Hence, using the Gaussian tail probability inequality with $|\tilde{f}(m)|\leq 1$, $|\psi_{j, k, m}|\leq 2^{-j/2}$, \fr{omeg} and \fr{Thres2}, as $\varepsilon, \delta \rightarrow 0$, simultaneously, yields
\be 
P_1\leq 2\phi\left(\frac{\gamma}{2} \lambda^{\alpha}_{\varepsilon, \delta}\left[\Var (\eta^{\varepsilon}_1)\right]^{-1/2}\right)= O\left(\varepsilon^{\frac{\gamma^2\rho_1^2}{16\sigma^2_{o1}}}\left[|\ln(\varepsilon)|\right]^{-1/2}\vee \delta^{\frac{\gamma^2\rho_2^2}{16\sigma^2_{o2}}}\left[|\ln(\delta)|\right]^{-1/2}\right),
\ee
and 
\be 
P_2\leq 2\phi\left(\frac{\gamma}{2} \lambda^{\alpha}_{\varepsilon, \delta}\left[\Var (\eta^{\delta}_2)\right]^{-1/2}\right) = O\left(\varepsilon^{\frac{\gamma^2\rho_1^2}{16\sigma^2_{o1}}}\left[|\ln(\varepsilon)|\right]^{-1/2}\vee \delta^{\frac{\gamma^2\rho_2^2}{16\sigma^2_{o2}}}\left[|\ln(\delta)|\right]^{-1/2}\right),
\ee
where $\phi(.)$ denotes the survival function of the standard normal random variable. This completes the proof of \fr{Largdev}. $\Box$\\
{\bf The proof of Theorem \ref{th:upperbds}.} The proof is very similar to the proof of Theorem 2 of Benhaddou~(2018a), and therefore we skip it. $\Box$\\

\end{document}